\documentclass{psp}

\usepackage{amsmath, amsfonts, amscd, amssymb}

\newcommand{\Z}{\mathbb{Z}}
\newcommand{\C}{\mathbb{C}}

\newcommand{\F}{\mathbb{F}}

\newcommand{\glnq}{GL(n,k)}
\newcommand{\glnk}{GL(n,k)}
\newcommand{\chg}{CH^*(BG)}

\newcommand{\MU}[1]{MU^*(#1)\hat{\otimes}_{MU^*}\Z}
\newcommand{\MUC}[2]{MU^*(#1)\hat{\otimes}_{MU^*} #2}
\newcommand{\MUL}[1]{MU^*(#1)\hat{\otimes}_{MU^*} \F_l}
\newcommand{\BPL}[1]{BP^*(#1)\hat{\otimes}_{BP^*} \F_l}

\newcommand{\tBPL}[1]{\td{BP}^*(#1)\hat{\otimes}_{BP^*} \F_l}
\newcommand{\td}[1]{\widetilde{#1}}
\newcommand{\gc}{G_\C}
\newcommand{\gq}{\Gamma_q}
\newcommand{\droite}{\longrightarrow}

\newnumbered{blank}{}[section]

\newunnumbered{note}{Note}

\newnumbered{rmk}[blank]{Remark}
\newnumbered{ex}[blank]{Example}
\newnumbered{notations}[blank]{Notations}

\title{Chow rings and cobordism of some Chevalley groups}
\author[Pierre Guillot]{PIERRE GUILLOT\\
Department of Pure Mathematics and Mathematical Statistics, \addressbreak
University of Cambridge, England\addressbreak
email\textup{: \texttt{p.guillot@dpmms.cam.ac.uk}}}

\receivedline{ Received \textup{23} September \textup{2002} }

\begin{document}

\maketitle

\begin{abstract}
We compute the cobordism rings of the classifying spaces of a certain class of Chevalley groups. In the particular case of the general linear group, we prove that it is isomorphic to the Chow ring.
\end{abstract}

\tableofcontents

\section{Introduction}
\begin{blank}\label{intro} If $G$ is an algebraic group over the complex numbers, then its classifying space $BG$ can be approximated by algebraic varieties, as was proved by Totaro in \cite{totaro}. Let alone the generalisation of the concept of classifying spaces to any field, this provides a number of new ways of studying $BG$, for example by investigating some algebro-geometric invariants. Among these, the Chow ring (the ring of algebraic cycles modulo rational equivalence) appears to be of some significance.

In most examples for which we can compute $CH^*(BG)$, the natural map to the ordinary cohomology ring tends to be injective, and $CH^*(BG)$ is often the ``nice'' part of $H^*(BG,\Z)$. This is most clearly illustrated by abelian groups. If $p$ is an odd prime, we have $H^*(\Z/p,\F_p)=\F_p[x,y]$ with $x$ of degree $2$ and $y$ of degree $1$, so that $\F_p[x]$ is just a polynomial ring, but $y$ on the other hand generates an exterior algebra, ie squares to $0$. Here $CH^*(B\Z/p)\otimes_\Z \F_p$ injects as the subring $\F_p[x]$. As we take products of such groups, the cohomology ring grows in complexity, but the Chow ring does not as much: $H^*((\Z/p)^n,\F_p)$ is the tensor product of a polynomial ring on $n$ variables and an exterior algebra on $n$ odd dimensional generators; the Chow ring is the polynomial part of that. This is not the ``even'' subring, but indeed in some sense the ``simple'' part of it. One of the goals of the present article is to show that something very similar holds in the case $G=GL(n,\F_q)$.

Another remarkable fact about the Chow ring is the factorization of the cycle map through the cobordism ring, that is, one has
$$\begin{CD}
CH^*(BG) @>cl^0>> \MU{BG} @>\pi >> H^*(BG,\Z)
\end{CD}$$ and $\pi\circ cl^0=cl$, the usual cycle map. Here $MU$ stands for complex cobordism. It is conjectured by Totaro that $cl^0$ is an isomorphism, and this is true in all known examples. When this is the case, and when the cycle map is injective, we have thus a purely topological definition as well as one from algebraic geometry of a certain subring of $H^*(BG,\Z)$ which we expect to be what we obtain by getting rid of the pathological (eg, nilpotent) elements. The mystery about this ring is the motivation to compute $CH^*(BG)$ or $\MU{BG}$ or both on new examples. It is also a first step towards the computation of the cohomology ring, which is very hard in general. \end{blank}

\begin{blank}
In this paper we investigate the case when $G$ is a Chevalley group, which means essentially a group of matrices over a finite field (precise definition in \ref{defchev}), such as $GL_n(\F_q)$, $SL_n(\F_q)$ or $Sp_{2n}(\F_q)$. We find that the cobordism ring of such a group, with mod $l$ coefficients, is for most choices of $l$ the same as that of the corresponding group over $\C$ (ie, $GL_n(\C)$, $SL_n(\C)$ or $Sp_{2n}(\C)$), provided that the finite field we work with contains the $l$-th roots of unity. Otherwise, some relations need to be introduced, corresponding to the action of a certain Galois group, and we give them explicitly. All this is the object of section \ref{section: cobordism}, see in particular theorem \ref{thm: cobordism} there.

In section \ref{section: chow} we specialise to the case of $GL_n$, for which we compute the Chow ring $CH^*(BGL_n(\F_q))$ (theorem \ref{thm: main}), and prove that it maps isomorphically to the cobordism ring (\ref{thm: iso}).

Our exposition starts with section \ref{section: prelim}, which summarizes a few important facts about the cohomology theories that we are using, and establishes a few lemmas of a rather technical nature.

\end{blank}

\begin{blank}
The general strategy that we have followed in our computations could probably be used to good effect in other situations. Namely, a decent starting point in many cases seems to have a look at the Morava K-theories first, which are very flexible and relatively easy to deal with, and then deduce some consequences on the Brown-Peterson cohomology via the very nice and powerful results obtained by Ravenel-Wilson-Yagita \cite{bpfrommorava}.
In this paper for example we have relied on the work of Tanabe \cite{moravachevalley} on Morava K-theories which we have adapted to $BP^*$. It is hopefully a good illustration of the method, showing how much can be easily done now that \cite{bpfrommorava} is available.

\end{blank}

\begin{blank}{\em Notations \& Conventions.} We will call $p$ and $l$ two distinct primes, and $q$ will be a power of $p$. The field with $q$ elements will be denoted by $\F_q$, and $\F$ stands for an algebraic closure of $\F_p$. For each field under discussion, $\mu_l$ will be the group of $l$-th roots of unity.

If $A$ is a ring with an action of a group $\Gamma$, then $A_{\Gamma}=A/(a - \gamma\cdot a)$ is the ring of coinvariants -- this is not quite standard, as the notation refers usually to the abelian group obtained by dividing out by the {\em subgroup} spanned by the $a - \gamma\cdot a$, as opposed to the {\em ideal} that they generate. But we choose to follow Tanabe \cite{moravachevalley}.

 On the other hand $A^{\Gamma}=\{a:\gamma\cdot a = a \}$ is classically the ring (resp. group if $A$ is only a group, etc...) of invariants.

Wreath products will be written on the right, eg $G\wr \Z/l$ or $G\wr S_n$.

A {\em space} will always mean a CW complex with finitely many cells in each dimension, although we will repeat this occasionally for emphasis. If $E$ and $F$ are spectra, then we put $E^*(F)=[F,E]$; if $X$ is a space, we put $E^*(X)=E^*(\Sigma^\infty (X^+))$; if $X$ is a pointed space, we put $\td{E}^*(X)=E^*(\Sigma^\infty (X))$. Here $\Sigma^\infty$ is the canonical functor from pointed spaces to spectra and $X^+$ is the disjoint union of $X$ and a base point. Putting $E^*=E^*(point)$, we have thus $E^*(X)=E^*\oplus \td{E}^*(X)$ and $\td{E}^*(X)=\ker (E^*(X) \to E^*)$.

If $X$ is a space and $E$ a spectrum, we will say that ``$X$ has no $\lim^1$ term with respect to $E$'' if $E^*(X)=\lim_\leftarrow E^*(X^m)$. Here and elsewhere the superscript refers to the $m$-th skeleton. This terminology is justified by Milnor's exact sequence (cf \cite{rudyak}, corollary 4.18).

\end{blank}

\section{Preliminaries on cohomology theories} \label{section: prelim}

\begin{blank}
We shall encounter several (generalised) cohomology theories: complex cobordism $MU^*$ and Brown-Peterson cohomology $BP^*$ (cf \cite{rudyak} chapter 7), Morava K-theories $K(j)^*$ (\cite{rudyak} chapter 9), and we will also mention $P(1)^*$ briefly in the course of one proof ({\em loc. cit.}). These spectra are to be taken at the prime $l$. So we have

$$\begin{array}{rcl}
MU^* & = & \Z[x_1,x_2,\cdots] \\
BP^* & = & \Z_{(l)}[v_1,v_2,\cdots] \\
P(1)^* & = & \F_l[v_1,v_2,\cdots] \\
K(j)^* & = & \F_l[v_j,v_j^{-1}]
\end{array}$$ with $x_i$ of degree $-2i$ and $v_i$ of degree $-2(l^i-1)$.

\end{blank}

\begin{blank}{\em Completed tensor products.}\label{complete}
For each such theory $h^*$ and spaces $X_i$, $i=1,2$ satisfying 
$$h^*(X_i)=\lim_{\leftarrow m} h^*(X^m_i)$$ where $X^m_i$ is the $m$-th skeleton of $X_i$, we will use the completed tensor products:

$$h^*(X_1)\hat{\otimes}_{h^*} h^*(X_2)=\lim_{\leftarrow m} [im(h^*(X_1)\to h^*(X^m_1)) \otimes_{h^*} im(h^*(X_2)\to h^*(X^m_2))]$$ 
and $h^*(X)\hat{\otimes}_{h^*} R$ is defined similarly, for a ring $R$ with a map $h^*\to R$ (by considering $R$ as graded but concentrated in dimension zero). Examples in view are the natural maps $MU^*\to \Z \to \Z_{(l)}\to \F_l$ and $BP^*\to\Z_{(l)} \to \F_l$. Recall that $$MU^*(X) \otimes_{MU^*} \Z_{(l)}=BP^*(X) \otimes_{BP^*} \Z_{(l)}$$ for, say, a finite dimensional CW complex $X$ so in particular
$$\MUL{X}=\BPL{X}$$ for a CW complex $X$.
\end{blank}

\begin{blank}
It is worth making a few comments about this last ring. Intuitively, it is obtained from $MU^*(X)$ (or $BP^*(X)$) by setting to zero all elements in the ideal generated by $l$ and the $v_i$'s, or rather in the closure of this ideal. For some nice spaces we can make this idea precise, namely, for spaces with no $\lim^1$ term with respect to $MU$ (or $BP$). For future reference we state this as a lemma:

\begin{lem*} \label{unwindtensor}
If $MU^*(X)=\lim MU^*(X^m)$, then the natural map $$MU^*(X)\to\MUL{X}$$ is surjective. It follows that $\MUL{X}$ is obtained from $MU^*(X)$ by dividing out by the closure of the ideal generated by the $v_i$'s and by $l$.

A similar result holds for $BP^*$. 
\end{lem*}

This follows from the arguments in \cite{totaro}, section 2. In fact something much stronger is proved there.

\end{blank}

\begin{blank}{\em An exactness result.} We shall need the following rather technical result later on. This is the only time when we shall use $P(1)^*$. It partially generalizes the half-exactness of the usual tensor product functor.

\begin{lem*} \label{exactness}
For $i=1,2,3$, let $X_i$ be a space with no $\lim^1$ term with respect to $BP$, and with even Morava K-theories. Suppose given an exact sequence (where all maps are induced by maps of tolopogical spaces):

$$\td{BP}^*(X_1)\to \td{BP}^*(X_2)\to \td{BP}^*(X_3) \to 0$$
Then there is an exact sequence

$$\tBPL{X_1}\to \tBPL{X_2}\to \tBPL{X_3} \to 0$$

\end{lem*}

\begin{proof}
Because $X_i$ has even Morava K-theories, the reduction mod $l$ of $BP^*(X)$ is $P(1)^*(X)$ (\cite{bpfrommorava}, theorem 1.9). Consider then the following diagram with exact columns:

$$\begin{CD}
    @.                0    @.            0      @.        \\
@.                  @VVV                @VVV         @.   \\
            @.       \ker_2  @>>>  \ker_3      @>>>  0    \\
@.                  @VVV                @VVV         @.   \\
\td{P(1)}^*(X_1) @>>> \td{P(1)}^*(X_2) @>>>  \td{P(1)}^*(X_3) @>>>  0    \\
@VVV                @VVV                @VVV         @.    \\
\tBPL{X_1}   @>>> \tBPL{X_2}   @>>>  \tBPL{X_3}   @>>>  0     \\
    @.              @VVV                @VVV         @.   \\
            @.       0       @.          0     @.
\end{CD}$$

By the lemma in \ref{unwindtensor}, we can define $\ker_i$ to obtain exact columns. We prove now that the rows are all exact.

For the middle one, simply apply $-\otimes_{\Z} \F_l$ to the given exact sequence.

From \ref{unwindtensor} again, we see that $\ker_i$ can be described as the closure of the ideal generated by the $v_i$'s in $\td{P(1)}^*(X_i)$. It follows that $\ker_2\to\ker_3$ has a dense image. But now, since $X_i$ is of finite type, each graded piece of $\td{P(1)}^*(X_i^m)$ is finite, so that $\td{P(1)}^*(X_i)=\lim_{\leftarrow} \td{P(1)}^*(X_i^m)$ (say, by Mittag-Leffler), and the natural topology on $\td{P(1)}^*(X_i)$ is compact Hausdorff. (This is the chief reason for using $P(1)^*$ instead of $BP^*$.) Thus $\ker_2\to\ker_3$ is surjective, and the first row is exact.

There remains only a trivial diagram chase to prove that the bottom row is exact, too.\end{proof}

\end{blank}

\begin{blank}{\em K\"unneth formulae.} \label{kunneth} A quick word on the cohomology of a product $X\times Y$ of two spaces. The situation for Morava K-theories is, as usual, very simple, and it is an understatement to say that the following result is well-known; however for lack of a good reference we submit a proof.

\begin{lem*} Let $X$, $Y$ be any two spaces. Then:
$$K(j)^*(X\times Y)=K(j)^*(X)\hat\otimes_{K(j)^*} K(j)^*(Y)$$
\end{lem*}

\begin{proof}(cf \cite{bpfrommorava}, proof of theorem 1.11.) Write $E$ for $K(j)$, $\otimes$ for $\otimes_{E^*}$, $\lim$ for inverse limits, and $E^*(X)_s$ for $im(E^*(X)\to E^*(X^s))$. We recall that $E^*$ is a graded field in the sense that every graded module over it is free.

Now if $Y$ is a finite complex we have $$E^*(X)\otimes E^*(Y)=E^*(X\times Y)$$ because if we see both sides of this equation as functors on the category of finite complexes, fixing $X$, then they are cohomology theories which agree on $S^0$ and with a natural transformation between them, so the result follows from \cite{rudyak}, proposition 3.19(i).

There is never a $\lim^1$ term for Morava K-theory (and for spaces), see \cite{bpfrommorava}, corollary 4.8. So in particular for any $X$ and $Y$ we have $\lim^1_i E^*((X\times Y)^i)=0$. From the naturality of the Milnor sequence applied to the subcomplexes $(X\times Y)^i$ and $X\times Y^i$ we draw $\lim^1_i E^*(X \times Y^i)=0$ and hence $$E^*(X\times Y)=\lim E^*(X\times Y^i)=\lim E^*(X)\otimes E^*(Y^i)$$ from the previous equation.

Observe that for a finitely generated $E^*$-module $M$ and any inverse system $\{ A_i \}$ of $E^*$-modules, we have $$(\lim A_i)\otimes M=\lim (A_i\otimes M)$$ because $M$ being free, this is just the statement that inverse limits commute with finite direct sums. We apply this below with $M=E^*(Y^i)$, $M=E^*(X)_s$ and $M=E^*(Y)_s$. The following computation should now be straightforward:

$$\begin{array}{rcl}
E^*(X\times Y) & = & \lim_j E^*(X)\otimes E^*(Y^j) \\
                & = & \lim_j (\lim_i E^*(X)_i) \otimes E^*(Y^j) \\
                & = & \lim_j \lim_i E^*(X)_i \otimes E^*(Y^j) \\
                & = & \lim_i \lim_j E^*(X)_i \otimes E^*(Y^j) \\
                & = & \lim_i E^*(X)_i \otimes (\lim_j E^*(Y^j)) \\
                & = & \lim_i E^*(X)_i \otimes (\lim_j E^*(Y)_j) \\
                & = & \lim_i \lim_j E^*(X)_i \otimes E^*(Y)_j \\
                & = & \lim_{i,j} E^*(X)_i \otimes E^*(Y)_j \\
                & = & E^*(X)\hat\otimes E^*(Y)
\end{array}
$$
\end{proof}

For $BP$ cohomology, there is the following nice result taken from \cite{bpfrommorava} (theorem 1.11): if $X_i$ is a CW complex of finite type, with even Morava K-theories, and such that $BP^*(X_i)=\lim BP^*(X^m_i)$ (for $i=1,2$) then
$$ BP^*(X\times Y) =  BP^*(X)\hat{\otimes}_{BP^*} BP^*(Y)$$
In particular, one has a map $$BP^*(X)\otimes_{BP^*} BP^*(Y)\to BP^*(X\times Y)$$ with a dense image, and similarly for Morava K-theory.

We will give some related results for Chow rings in \ref{chowkunneth}.

\end{blank}

\begin{blank}{\em More exactness.}\label{exactrings} To finish these preliminaries, we have the following proposition dealing with a certain type of exact sequences. Its formulation might seem cumbersome at first; the reason is that we want to exploit the ring structure on our cohomology groups, and the notion of an exact sequence of rings is delicate to phrase (the image of a homomorphism is always a subring, which is almost never an ideal at the same time). Here we can take advantage of the decomposition $E^*(X)=E^* \oplus \tilde{E}^*(X)$ using reduced cohomology and make some sense (however we shall not use the terminology of ``augmented rings'').

\begin{prop*} Suppose given spaces $X$, $Y$ and $Z$ with no $\lim^1$ term with respect to $BP$, and suppose that $X\times Y$ has no $\lim^1$ term either. Assume given maps 
$$\begin{CD}
X @<f<< Y @<g<< Z
\end{CD}$$ with $f\circ g$ inessential, giving rise to an exact sequence

$$\begin{CD}
\td{K(j)}^*(X)\otimes_{K(j)^*} K(j)^*(Y) @>f^*\otimes id>> \td{K(j)}^*(Y) @>g^*>> \td{K(j)}^*(Z) @>>> 0 
\end{CD}$$ for all $j>0$. Then the same is true with $K(j)^*(-)$ replaced by $\BPL{-}$.

\end{prop*}

\begin{rmk*} The exact sequence says that the kernel of $g^*$ is the {\em ideal} generated by the image of $f^*$. Thus the multiplicative structure plays the prominent role here.
\end{rmk*}

\begin{proof}
Write $E$ for $K(j)$. We use the Kunneth formula for Morava K-theory:
$$\begin{array}{rcl}
E^*(X \times Y) & = & E^*(X)\hat\otimes_{E^*} E^*(Y) \\
                & = & E^* \oplus \td{E}^*(Y)  \oplus \td{E}^*(X)\hat\otimes_{E^*} E^*(Y)
\end{array}$$
Note that the sum of the last two terms above is $\td{E}^*(X\times Y)$, and the  last summand is the kernel of the map $\td{E}^*(X\times Y)\to\td{E}^*(Y)$ induced by the inclusion $Y\to X\times Y$. Let $X\propto Y$ be the cone of this map (ie $X\times Y/*\times Y$). The cofibration gives an exact sequence on reduced cohomology, and the projection $X\times Y \to Y$ gives a splitting, so that 
$$\begin{array}{rcl}
\td{E}^*(X\propto Y) & = & \ker \td{E}^*(X\times Y)\to\td{E}^*(Y) \\
                     & = & \td{E}^*(X)\hat\otimes_{E^*} E^*(Y)
\end{array}$$
We conclude that the map $F: Y\to X\propto Y$ induced by $(f,id)$ gives rise to an exact sequence  
$$\begin{CD}
\td{K(j)}^*(X \propto Y) @>>> \td{K(j)}^*(Y) @>>> \td{K(j)}^*(Z) @>>> 0 
\end{CD}$$ for all $j>0$. Also, $F\circ g$ is inessential. Therefore by theorem 1.18 in \cite{bpfrommorava}, we have also
$$\begin{CD}
\td{BP}^*(X \propto Y) @>>> \td{BP}^*(Y) @>>> \td{BP}^*(Z) @>>> 0 
\end{CD}$$ and by (\ref{exactness}, lemma) we have
$$\begin{CD}
\tBPL{X \propto Y} @>>> \tBPL{Y} @>>> \tBPL{Z} @>>> 0 
\end{CD}$$
Now there are maps
$$BP^*(X)\otimes_{BP^*} BP^*(Y) \to BP^*(X\times Y) \to \BPL{X\times Y}$$
The first map has a dense image by (\ref{kunneth}), and the second map is surjective by (\ref{unwindtensor}). Arguing as above we end up with a map $$\td{BP}^*(X)\otimes_{BP^*} BP^*(Y) \to \tBPL{X\propto Y}$$ with a dense image. Clearly we also have a map $$\left[ \tBPL{X}\right] \times \left[ \BPL{Y} \right] \to \tBPL{X\propto Y}$$ with a dense image, but here by compactness the map is actually surjective. 

Finally, we patch up things together and obtain the desired exact sequence
$$\begin{CD}
\left[ \tBPL{X}\right] \otimes_{\F_l} \left[ \BPL{Y}\right] @>f^*\otimes id>> \tBPL{Y} @>>> \tBPL{Z} @>>> 0 
\end{CD}$$
This completes the proof.
\end{proof}

\begin{rmk*} During the course of this proof we have implicitly used the fact that the space $X\propto Y$ has no $\lim^1$ term for $BP$. This follows easily from the injectivity of $$\td{BP}^*(X\propto Y) \to \td{BP}^*(X\times Y)$$ and the naturality of Milnor's exact sequence.
\end{rmk*}

\end{blank}

\section{Cobordism of Chevalley groups} \label{section: cobordism}

\begin{blank} \label{defchev}
Let $G$ be a connected, reductive, and split group scheme over $\Z$. Here split means that there exists a maximal torus $T$ of $G$ which is {\em defined} and {\em split} (ie isomorphic to a product of copies of $G_m$) over $\Z$. We are interested in the finite groups $G(\F_q)$, usually called {\em Chevalley groups}.

 Then $\gc=G\times_{Spec(\Z)} Spec( \C)$ is viewed as a Lie group (and in fact the notation will occasionally refer to a maximal compact Lie subgroup, which should not cause any confusion as the inclusion map induces a homotopy equivalence; the two classifying spaces also share the same homotopy type). 

 We shall work under the assumption that $H^*(\gc,\Z)$ has no $l$-torsion.

\end{blank}

\begin{blank}{\em The associated Lie group.} \label{blank: lie}
The study of $\gc$ is quite easy, given the well-known results obtained by Borel for ordinary cohomology. We have:

\begin{prop*}
Let $n$ be the rank of $G$. There are elements $s_i$, for $1\le i \le n$, of positive even degree, such that:
$$BP^*(B\gc)=BP^*[[s_1,\cdots,s_n]]$$
It follows that $$\BPL{B\gc}=H^*(B\gc,\F_l)=\F_l[s_1,\cdots,s_n]$$
 
\end{prop*}

\begin{proof}
This follows from the classical result (cf \cite{borelphd}) $$H^*(B\gc,\Z_{(l)})=\Z_{(l)}[s_1,\cdots,s_n]$$ and similarly with $\F_l$ in place of $\Z_{(l)}$, by using the Atiyah-Hirzebruch spectral sequence.
\end{proof}
\end{blank}

\begin{blank}{\em Brauer lifts.} \label{blank: brauer}
In order to study now the finite groups $G(\F_q)$ we shall make use of ``Brauer lifts'', a general term meaning that we lift some representations of a group from characteristic $p$ to characteristic $0$. A possible approach is described in the proposition below; later in \ref{morebrauer} another point of view will be preferable.

\begin{prop*}
There exists a map $$BG(\F)\droite B\gc$$ natural with respect to maps of group schemes, such that the maps $$BP^*(B\gc)\droite BP^*(BG(\F_q))$$ are all surjective, for all choices of $q$. It follows that $$\BPL{B\gc}\droite \BPL{BG(\F_q)}$$ is also surjective. 

\end{prop*}

\begin{proof}
The map is described in \cite{brauercohom}. In \cite{moravachevalley}, it is shown that it induces surjective maps $$K(j)^*(B\gc)\droite K(j)^*(BG(\F_q))$$ for all $j>0$ (and in this case the kernel is also known). By a result of Ravenel-Wilson-Yagita \cite{bpfrommorava}, this implies that the induced map on $BP$ cohomology is surjective as well. The final statement follows by the right exactness of the (completed) tensor product.
\end{proof}

\begin{rmk*}[1]
Note that for classifying spaces of compact Lie groups, such as $B\gc$ and $BG(\F_q)$, there is no $\lim^1$ term for $BP$ cohomology (and for the skeletal filtration). This is essential to use the results of \cite{bpfrommorava} in their strong form (ie, with $BP$ itself instead of its completion). This is a major concern which explains why very few of our statements will involve $BG(\F)$, for which there might be {\em a priori} a non-vanishing $\lim^1$.

\end{rmk*}

\begin{rmk*}[2] The kernel of the map above can in fact be described, but we postpone this to \ref{galois}, in order to prove as much as possible by elementary means.
\end{rmk*}

\end{blank}

\begin{blank}{\em Tori.}
Let $K$ be a finite field of characteristic $p$ which contains the $l$-th roots of unity and let $T$ be a torus (considering our framework (cf \ref{defchev}), this means that $T$ is split over $\Z$ and hence over $K$). Then $T(K)$ is a product of cyclic groups whose order is a multiple of $l$. Therefore
$$\BPL{BT(K)}=\F_l[\eta_1,\cdots,\eta_n]$$ where $n$ is the dimension of $T$, and where the $\eta_i$'s have degree $2$. But of course $\BPL{BT_\C}$ has a similar description (these are all well-known results). Consequently the surjective map between this two rings that we obtained in (\ref{blank: brauer}, proposition) is in fact an isomorphism.

\end{blank}

\begin{blank} \label{injecttorus}
Quite fortunately, we can relate easily any group $G$ as above to its maximal torus, in terms of cohomology. More precisely:

\begin{prop*}
The restriction map
$$\BPL{B\gc}\droite \BPL{BT_\C}$$
is injective. Moreover if the order of the Weyl group is prime to $l$, then the image of this map is precisely the subring of invariants.
\end{prop*}

\begin{proof}
As observed above (\ref{blank: lie}, proposition), this is simply the map
$$H^*(B\gc,\F_l)\droite H^*(BT_\C,\F_l)$$ It is well known, since $BG$ has no $l$-torsion, that the restriction map gives an isomorphism $$H^*(B\gc,\Z_{(l)})\droite H^*(BT_\C,\Z_{(l)})^W$$ where $W$ denotes the Weyl group. Consider the exact sequence of coefficients:
$$\begin{CD}
0 @>>> \Z_{(l)} @>\times l>> \Z_{(l)} @>>> \F_l @>>> 0  
\end{CD}$$

This gives rise to exact sequences in cohomology:
$$\begin{CD}
0 @>>> H^k(B\gc,\Z_{(l)}) @>\times l>> H^k(B\gc,\Z_{(l)}) @>>> H^k(B\gc,\F_l) @>>> 0  
\end{CD}$$ using either the fact that the rings are concentrated in even dimensions, or that multiplication by $l$ is clearly injective here. Apply this to $T$ as well, appeal to naturality, and apply the functor $(-)^W$ to obtain the following commutative diagram:

{\small
$$\begin{CD}
0 @>>> H^k(B\gc,\Z_{(l)}) @>>> H^k(B\gc,\Z_{(l)}) @>>> H^k(B\gc,\F_l) @>>> 0  \\
@. @VV\approx V @VV\approx V @VVV \\
0 @>>> H^k(BT_\C,\Z_{(l)})^W @>>> H^k(BT_\C,\Z_{(l)})^W @>>> H^k(BT_\C,\F_l)^W @>>> H^1(W, H^k(BT_\C,\Z_{(l)}))  
\end{CD}$$
}

Injectivity is now proved by a quick diagram chase. If the order of $W$ is prime to $l$, then $H^1(W, H^k(BT_\C,\Z_{(l)}))=0$ and the result follows.
\end{proof}

\end{blank}

\begin{blank}{\em Big finite fields.} \label{blank: bigones} We can now finish off the computations for $BG(K)$ if $K$, as above, is a finite field of char $=p$ containing the $l$-th roots of unity. Indeed, putting together the information of the previous paragraphs, we get a commutative and exact diagram:
$$\begin{CD}
  @.              @.       0 \\
@.         @.             @VVV  \\
0 @>>> \BPL{B\gc} @>>> \BPL{BT_\C} \\
@.         @VVV           @VVV  \\
  @.   \BPL{BG(K)} @>>> \BPL{BT(K)} \\
@.         @VVV           @VVV   \\
  @.        0      @.         0
\end{CD}$$

From this it follows that {\em for $G$ as in the introduction and $K$ as above, the map $$\BPL{B\gc}\droite \BPL{BG(K)}$$ is an isomorphism.}

\end{blank}

\begin{blank}{\em Galois actions.} \label{galois}
If now $k=\F_q$ is any field, let $K=k(\mu_l)$ (extension obtained by adding the $l$-th roots of unity) and let $\Gamma_q=Gal(K/k)$, which is cyclic generated by the Frobenius automorphism $\gamma$. Let $n$ be the rank of $G$ and let $r=[K:k]$.

We have a surjective map $$\BPL{BG(K)} \droite \BPL{BG(k)}$$ from (\ref{blank: brauer}, proposition) and (\ref{blank: bigones}). Moreover the Galois group $\gq$ acts on both groups, trivially on $\BPL{BG(k)}$. Therefore there is actually a surjective map from $(\BPL{BG(K)})_{\gq}$ (coinvariants) to $ \BPL{BG(k)}$ and one has:

\begin{prop*} This map is an isomorphism, ie 
$$(\BPL{BG(K)})_{\gq}=\BPL{BG(k)}$$
\end{prop*}

\begin{proof}
Tanabe proved that this was true for Morava K-theory in \cite{moravachevalley}, which can be stated by saying that we have an exact sequence just as in (\ref{exactrings}, proposition) with $X=Y=BG(K)$, $Z=BG(k)$, and $f^*(a)=a - \gamma\cdot a$, keeping the same notations. This very proposition will give us the result we want if we can only prove that we can find a map $f$ of topological spaces inducing such an $f^*$. But this is easy: the maps $a\mapsto a$ and $a\mapsto \gamma\cdot a$ being both induced by maps of topological spaces, so is their difference, at least after suspending which has no effect on cohomology\footnote{\scriptsize For spectra $E$ and $F$, $[E,F]$ is always an abelian group, so that given two spaces $X$ and $Y$ one is tempted to consider $E=\Sigma^\infty(X^+)$ and $F=\Sigma^\infty(Y^+)$ if one wants to take the difference between two maps $X\to Y$. However one must be careful, as most results in \cite{bpfrommorava}, on which we rely heavily, are strictly unstable, ie do not apply to spectra and their maps.}.
\end{proof}

\end{blank}

\begin{blank} So let us describe the ring of coinvariants a bit more precisely. Starting with an $n$-dimensional torus as usual, we have
$$\BPL{BT(K)}=\F_l[\eta_1,\cdots,\eta_n]$$ and the description of the $\eta_i$'s as Euler classes makes it clear that $\gamma$ sends $\eta_i$ to $q\cdot \eta_i$. 

For general $G$, the ring $$\BPL{BG(K)}=\F_l[s_1,\cdots,s_n]$$ injects into the analogous ring for a maximal torus, and $s_i$ is sent to a homogenous polynomial in the $\eta_j$'s of degree $|s_i|/2$; therefore $\gamma$ sends $s_i$ to $q^{|s_i|/2}\cdot s_i$. We have then the following lemma (which will be used later on with $R^k=2k$ -th graded piece of $\BPL{BG(K)}$):

\begin{lem*}
 Suppose $R=\F_l[s_1,\cdots,s_n]$ is a graded polynomial ring with an action of a cyclic group $\Gamma=<\gamma >$ given by $\gamma\cdot s_i=q^{|s_i|}s_i$. Then the ring of coinvariants $R'=R/(a-\gamma\cdot a)$ is isomorphic to $\F_l[s_i: r \textrm{ divides } |s_i|]$.
\end{lem*}

\begin{proof} To see this, observe first that the condition $r|k$ is equivalent to the $l$-th roots of unity being in $\F_{p^k}$, and so in turn this is equivalent to $q^k=1$ mod $l$. Suppose now that $k=|s_i|$ is {\em not} a multiple of $r$, then $1-q^k$ is invertible in $\F_l$, and we can write
$$s_i=s_i\cdot\frac{1-q^k}{1-q^k}=\frac{s_i}{1-q^k} - \gamma\cdot\frac{s_i}{1-q^k}$$ so that $s_i$ maps to $0$ in $R'$. Hence there is a surjective map $\F_l[s_i: r \textrm{ divides } |s_i|]\to R'$. We leave to the reader the easy task of proving the injectivity of this map.
\end{proof}

\end{blank}

\begin{blank} \label{blank: muchev} To summarize, we have proved the following (casually replacing $BP$ with $MU$, cf (\ref{complete})):

\begin{thm*} \label{thm: cobordism}
Let $G$ be a reductive, connected group scheme (over $\Z$). Let $l$ be a prime number such that $H^*(\gc,\Z)$ has no $l$-torsion. Let $K$ be a finite field of characteristic $p\ne l$ containing the $l$-th roots of unity, and let $n$ be the rank of $G$. Then there are elements $s_i$ of even degree such that 
$$\begin{array}{rcl}
\MUC{BG(K)}{\F_l}  & = &  \MUC{B\gc}{\F_l}\\
                   & = & H^*(B\gc,\F_l) \\
                   & = & \F_l[s_1,\cdots,s_n]
\end{array}$$

 If $k$ is any finite field, let $K=k(\mu_l)$ and let $r=[K:k]$. Then 
$$\MUC{BG(k)}{\F_l}=\F_l[s_i:2r \textrm{ divides } |s_i|]$$

\end{thm*}

\end{blank}

\begin{blank}{\em Example: the symplectic group.}
Let us have a look at $G=Sp_{2n}$. We have
$$H^*(Sp_{2n}(\C),\Z)=\Lambda (e_1,e_2,\cdots,e_n)$$ and
$$H^*(BSp_{2n}(\C),\Z)=\Z[q_1,q_2,\cdots,q_n]$$ where $|e_i|=4i-1$ and $|q_i|=4i$. Say for example that we choose $l=5$ and that we look at fields of characteristic $3$. If we take $q=81$, the $5$-th roots of unity are already in $\F_{81}$, and we have
$$\MUC{BSp_{2n}(\F_{81})}{\F_5}=\F_5[q_1,\cdots,q_n]$$
If we take $q=9$ we need to go to $\F_{81}$ and $r=2$, but all degrees appearing above are divisible by $2r=4$. So we also obtain that $$\MUC{BSp_{2n}(\F_9)}{\F_5}=\F_5[q_1,\cdots,q_n]$$

However for $q=3$ we have $r=4$, and also for $q=27$ (the $5$-th roots of unity appear only in $\F_{531 441=27^4}$) and we drop the ``odd'' generators, whose degree is not divisible by $2r=8$. Hence $$\MUC{BSp_{2n}(\F_3)}{\F_5}=\MUC{BSp_{2n}(\F_{27})}{\F_5}=\F_5[q_2,q_4,\cdots]$$

\end{blank}

\section{The Chow ring of the general linear group over a finite field}
\label{section: chow}

\begin{blank} In order to compute now $CH^*(B\glnq)/l$, we shall follow closely Quillen's paper \cite{quillengln} where the cohomology ring is investigated. The argument is greatly simplified in the case of Chow rings thanks to the nicer K\"unneth formula at our disposal.

Throughout, $k$ will be a finite field of characteristic $p$ and $K$ will denote, as above, the extension of $k$ obtained by adding the $l$-th roots of unity. Again put $r=[K:k]$, and write the long division $n=rm + e$ with $0\le e < r$. We shall assume that $m\ge 1$, noting that in the trivial case when $m=0$ (ie, $n<r$), then the order of $GL(n,k)$ is prime to $l$: to see this, recall first that this order is 
$$ q^{n(n-1)/2}\prod_{i=1}^n (q^i-1)$$ and note that $r$ is the smallest integer such that $l$ divides $q^r-1$.

Hence in this case $CH^*(B\glnk)/l=\Z/l$ in dimension $0$ (and $0$ in other dimensions).

A last notation: write $q^r-1=l^ah$, with $h$ prime to $l$. By choice of $r$, we have $a\ge1$.
\end{blank}

\begin{blank}{\em A subgroup of $\glnk$.}
We let $C=K^\times$, a cyclic group of order $q^r-1$. We have an obvious representation of $C$ on $K$ which is $1$-dimensional. If we let $\pi=Gal(K/k)$, then we can extend the action of $C$ to one of the semi-direct product $C\rtimes\pi$; the representation we obtain is faithful. Taking direct sums, we now have a faithful representation of $(C\rtimes\pi)^m$ on $K^m$. 

In turn, we can form the semi-direct product with $S_m$ to obtain the wreath product $(C\rtimes\pi)\wr S_m=(C\rtimes\pi)^m\rtimes S_m$, and we can again extend the above representation to this new group. It is still faithful.

Regarding $K$ as an $r$-dimensional vector space over $k$, and adding an $e$-dimensional trivial representation, we end up with an embedding of $(C\rtimes\pi)^m\rtimes S_m$ in $\glnk$. It has the advantage of providing an embedding of $C^m$ in $\glnk$ such that the actions of $\pi^m$ and $S_m$ are realized by inner automorphisms of $\glnk$.

The subgroup $C^m$ enjoys the following property:

\begin{lem*}   \label{lem: allconj}
Any abelian subgroup of $\glnk$ of exponent dividing $l^a$ is conjugate to a subgroup of $C^m$.
\end{lem*}

A proof can be found in \cite{quillengln}, lemma 12.
\end{blank}

\begin{blank}{\em K\"unneth considerations.} \label{chowkunneth} Let $G$ be a finite group. We shall call (*) the following
condition: $BG$ can be cut into open subspaces of affine spaces,
and the map $\chg\to H^*(BG,\Z)$ is split injective.

\begin{lem*}[1] \label{lem: totaro}
Let $G$ be a finite group satisfying (*). Then the map
$$CH^*(B(G\wr  \Z/l))/l \to CH^*(BG^l)/l \oplus CH^*(B(G\times \Z/l))/l$$
is injective. Moreover, $G\wr  \Z/l$ satisfies (*).

\end{lem*}

This lemma is essentially taken from \cite{totaro}; the result is stated in this form at the end of section 9 for a special case, but it is clear from the material contained in sections 8 and 9 that it holds as we have stated it. See lemma 8.1 in particular for the last statement.
The following is also in loc. cit. in just a slightly different form:

\begin{lem*}[2] \label{lem: totarosoupedup}
Let $X$ and $Y$ be varieties over the complex numbers. We assume that $CH_*(X)\to H_*(X,\Z)$ is split injective, and that $CH_*(Y)\otimes \Z/l\to H_*(X,\Z)\otimes\Z/l$ is injective.
Then $$CH_*(X)\otimes CH_*(Y)\otimes \Z/l \longrightarrow CH_*(X\times Y)\otimes \Z/l$$ is injective. If moreover $X$ or $Y$ can be partitioned into open subsets of affine spaces, then $$CH_*(X)\otimes CH_*(Y)\otimes \Z/l \approx CH_*(X\times Y)\otimes \Z/l$$

\end{lem*}

\end{blank}

\begin{proof}
Note first that the map $CH_*(Y)\otimes \Z/l\to H_*(X,\Z)\otimes\Z/l$ is in fact {\em split} injective, it being a linear map of $\F_l$-vector spaces.
Therefore the map $$CH_*(X)\otimes CH_*(Y)\otimes \Z/l \longrightarrow H_*(X,\Z)\otimes H_*(Y,\Z)\otimes \Z/l$$ is injective. Also, because $$ H_*(X,\Z)\otimes H_*(Y,\Z)\longrightarrow H_*(X\times Y,\Z)$$ is split injective, the map
$$ H_*(X,\Z)\otimes H_*(Y,\Z)\otimes \Z/l \longrightarrow H_*(X\times Y,\Z)\otimes \Z/l$$ is injective.

The first assertion follows then from the commutative diagram:

$$\begin{array}{ccc}
CH_*(X)\otimes CH_*(Y)\otimes \Z/l & \longrightarrow & CH_*(X\times Y)\otimes \Z/l \\
\downarrow & & \downarrow \\
H_*(X,\Z)\otimes H_*(Y,\Z)\otimes \Z/l & \longrightarrow &  H_*(X\times Y,\Z)\otimes \Z/l 
\end{array}$$

The second part follows trivially from the fact that $-\otimes \Z/l$ is right exact (the map being surjective before tensoring, cf the arguments in \cite{totaro}).\end{proof}

\begin{blank}{\em Detecting families.}
A family of subgroups $(G_i)_{i\in I}$ of $G$ will be called an
{\em $l$-detecting family} if the map $$\chg/l \to \bigoplus_{i\in
I} CH^*(BG_i)/l$$ is injective. Throughout, the indexing sets $I$ will be finite.

\begin{lem*} \label{lem: detect}
Suppose that $G$ satisfies (*), and let $l$ be a prime number. If
$G$ has an $l$-detecting family of abelian subgroups whose
exponent divides $l^a$, then $G\wr S_n$ has the same
property.
\end{lem*}

\end{blank}

\begin{proof}(cf \cite{thomas}, corollary of 9.4)
We proceed by induction, the result being trivial for $n=1$. We
have two cases:

\begin{enumerate}

\item{$n=n'l$.} In this case the subgroup $(G\wr \Z/l)\int
S_{n'}$ has index prime to $l$, so it suffices to prove the
existence of an $l$-detecting family for this subgroup. By the
first part of (\ref{lem: totaro}, lemma (1)), $G\wr \Z/l$ has such a
family; by the second part of the same lemma, it satisfies (*), so
that we may simply appeal to the induction
hypothesis.

\item{$n$ is prime to $l$.} The subgroup $G\times (G\wr S_{n-1})$
has index prime to $l$, so it detects the mod $l$ Chow rings, and it is enough to exhibit an $l$-detecting family for this subgroup. Now by the induction hypothesis there is such a family $(H_i)_{i\in I}$ for $H=G\wr S_{n-1}$. Consider the following commutative diagram

$$\begin{array}{ccc}
CH^*(BH)\otimes \Z/l & \longrightarrow & \bigoplus_i CH^*(BH_i)\otimes \Z/l \\
\downarrow & & \downarrow \\
H^*(H,\Z)\otimes \Z/l & \longrightarrow & \bigoplus_i H^*(H_i,\Z)\otimes \Z/l
\end{array}$$
The top map is injective, as is the right vertical one (because $H_i$ is abelian). We deduce that the left vertical map is injective as well.

From (\ref{lem: totarosoupedup}, lemma (2)), we conclude that
$$ CH^*(B(G\times H))/l \approx CH^*(BG)\otimes CH^*(BH) \otimes \Z/l$$
and the same holds with $H_i$ in place of $H$. It follows that the map
$$ CH^*(B(G\times H))/l \longrightarrow \bigoplus_i CH^*(B(G\times H_i))/l$$
is none other than the map on the top of the last diagram above (which is split injective by arguments similar to the ones in the previous proof), tensored by $CH^*(BG)$. It is therefore injective. 

Let $(K_j)_{j\in J}$ be an $l$-detecting family for $G$, and repeat the argument with some $H_i$ in place of $G$, and the $K_j$'s in place of the $H_i$'s: we end up with an injective map:
$$ CH^*(B(G\times H_i))/l \longrightarrow \bigoplus_j CH^*(B(K_j\times H_i))/l$$
Composing, we finally obtain an injective map
$$ CH^*(B(G\times H))/l \longrightarrow \bigoplus_{i,j} CH^*(B(K_j\times H_i))/l$$
Thus $(K_j\times H_i)_{(i,j)\in I\times J}$ is an $l$-detecting family of abelian groups for $G\times H$.
\end{enumerate}

This completes the induction step.\end{proof}

\begin{blank}{\em More Brauer lifts.} \label{morebrauer}
We recall another way of lifting representations from characteristic $p$ to characteristic $0$, thus creating interesting new characters. From chapter 18 of Serre's book \cite{serrerep} we keep the following ingredient:

\begin{prop*}
Let $G$ be a finite group, and let $\rho : G\to GL(n,\F)$ be a representation of $G$ over a finite field $\F$. Let $m$ be the least common multiple of the orders of the regular elements of $G$, let $\mathbb{E}$ be the extension of $\F$ obtained by adjoining the $m$-th roots of unity, and suppose given an embedding $\psi : \mathbb{E}^* \to \C^*$. Put $$\chi (g)=\sum_{\lambda} \psi(\lambda)$$ where the sum runs over all eigenvalues of $\rho(g)$ counted with multiplicities. 

Then $\chi\in R(G)$, that is, $\chi$ is a virtual character of $G$. 
\end{prop*}

\begin{rmk*}
Recall that the {\em regular} elements of $G$ are those whose order is prime to $p$. Note also that such an embedding $\psi$ clearly always exists, but is not unique. Most of what follows depends on this choice. By abuse of language we refer to $\chi$ simply as the Brauer lift of $\rho$ (thus ignoring the dependance on $\psi$). 
\end{rmk*}

\end{blank}

\begin{blank} In the case of $G=\glnk$ there is an obvious candidate for such a lift to characteristic $0$, and it turns out that the Chern classes of the virtual character thus obtained satisfy an important property:

\begin{prop*} \label{prop: restrictchern} Write $CH^*(BC^m)/l=\Z/l[\eta_1,\cdots,\eta_m]$. Let $\rho$ be the lift, as defined above, of the obvious representation of $\glnk$ on $k^n$, and let $\bar{c_i}$ be the restriction to $C^m$ of the Chern class $c_i(\rho)$. Then $\bar{c_{ir}}=0$ if $i>m$, and otherwise $\bar{c_{ir}}$ is (up to sign) the $i$-th symmetric function on the variables $\eta_j^r$, for $1\le j \le m$.
\end{prop*}

\end{blank}

\begin{proof}
Consider first the cyclic group $C$ and its natural representation $\theta$ on $\C$. Its first Chern class $\eta=c_1(\theta)$ is such that $CH^*(BC)/l=\Z/l[\eta]$. Consider also the representation $$W=\theta\oplus\theta^{\otimes q}\cdots \oplus\theta^{\otimes q^{r-1}}$$ The homomorphism $x\mapsto qx$ of $C$ sends $\theta$ to $\theta^{\otimes q}$ and hence it sends $\eta$ to $q\eta$. The ring of fixed elements in the Chow ring is $\Z/l[\eta^r]$ (see the argument in the proof of theorem \ref{thm: main} below). Because $W$ is unchanged by this transformation, its Chern classes are necessarily $0$ in dimensions other than $r$ (and $0$). In dimension $r$, by the sum formula, $c_r(W)=q^{r(r-1)/2}\eta^r$. But note that $q^{r(r-1)/2}=(-1)^{r-1}$ modulo $l$.

Next consider the group $C^m$ whose Chow ring is $\Z/l[\eta_1,\cdots,\eta_m]$. We let $W_i$ be the lift of the above representation of $C$ via the $i$-th projection, and we now use the letter $W$ to denote $W=\bigoplus_i W_i$. By the formula giving the Chern class of a sum, we deduce that $c_{ir}$ is the $i$-th symmetric function in the $\eta_j^r$ (and is $0$ if $j>m$), up to a sign.

The proof will therefore be complete if we show that the restriction of $\rho$ to $C^m$ is isomorphic to $W$. To see this, let $x$ be a generator of $C$. The action of $x$ on $K$, considered as a $k$-linear map, has by definition the trace $$tr_{K/k}(x)=\sum_{\sigma\in Gal(K/k)} \sigma (x) = \sum_{i=1}^{r-1} x^{q^i}$$ Result follows.\end{proof}

\begin{blank}{\em Main result.} We are now in position to prove:

\begin{thm*} \label{thm: main}
Let $k=\F_q$ where $q$ is a power of the prime $p$. Let $l$ be an odd prime different from $p$, put $r=[k(\mu_l):k]$, and put $m=[n/r]$. Let $\rho$ be the Brauer lift of the natural action of $\glnk$ on $k^n$, and let $c_i$ be its $i$-th Chern class. Then:

$$CH^*(B\glnk)/l=\Z/l[c_r,c_{2r},\cdots,c_{mr}]$$

\end{thm*}

Note that if $m=0$, this means $CH^*(B\glnk)/l=\Z/l$ in dimension $0$, as announced in the introduction.

\end{blank}

\begin{proof}
Let $N=\pi\wr  S_m$, a subgroup of $\glnk$ which normalizes $C^m$. We have a restriction map:
$$CH^*(B\glnk)/l\longrightarrow (CH^*(BC^m)/l)^N$$
which we will prove to be an isomorphism.

Let us write $CH^*(BC^m)/l=\Z/l[\eta_1,\cdots,\eta_m]$. Recall that $\pi$ is cyclic generated by the Frobenius map $x\mapsto x^q$; this element acts on $\eta_i$ by multiplication by $q$, and hence on $\eta_i^s$ by multiplication by $q^s$. It follows that an invariant element of the Chow ring has to be a polynomial in the $\eta_i^r$ (because $q^s=1$ modulo $l$ if and only if the $l$-th roots of unity are in $\F_{q^s}$, ie if and only if $r|s$). Furthermore $S_m$ acts by permuting the indices, so that finally the ring of fixed elements is the polynomial algebra on the $s_i$, the symmetric functions in the $\eta_j^r$.

It is then immediate from (\ref{prop: restrictchern}, proposition) that the above map is surjective. If we can show that it is injective as well, then $CH^*(B\glnk)/l$ will have the desired description.

First consider the group $C\wr  S_m$. Now, $C$ has a cyclic subgroup of order $l^a$, which is an $l$-detecting family in itself, so that by (\ref{lem: detect}, lemma), $C\wr  S_m$ has also an $l$-detecting family of abelian groups whose exponent divides $l^a$ (noting that a cyclic group satisfies (*)). 
Observe now that the index of $C\wr  S_m$ in $\glnk$ is prime to $l$ (cf proof of proposition 4 in \cite{quillengln}, so it detects the modulo $l$ Chow ring of $\glnk$, and therefore the family of subgroups just considered works for $\glnk$ itself.
But now from (\ref{lem: allconj}, lemma), it is clear that the restriction map above is injective.\end{proof}

\begin{coro*}
If $k$ contains the $l^b$-th roots of unity for some integer $b$, then
$$CH^*(B\glnk)/l^b=\Z/l^b[c_1,\cdots,c_n]$$
\end{coro*}

\begin{proof}
For $b=1$ this is just a particular case of the above theorem. For general $b$, it is always true that $CH^*(BC^n)/l^b=\Z/l^b[\eta_1,\cdots,\eta_n]$. Since the Galois group acts now trivially, the ring of fixed elements remains easy to describe: it is $\Z/l^b[s_1,\cdots,s_n]$ (notation as above). Therefore the restriction map already considered is still surjective even if we take the coefficients in $\Z/l^b$. One concludes using induction and a cardinality argument.
\end{proof}

\begin{blank}{\em Back to cobordism.}
Since we have seen that $CH^*(B\glnk)/l$ and $\MUL{B\glnk}$ have the same description, it is quite easy to see (at least) that the cycle map is an isomorphism: indeed consider the diagram:

$$\begin{CD}
CH^*(BGL_n(\C))/l  @>>> CH^*(B\glnk)/l \\
@V\approx VV            @VVcl^0 V \\
\MUL{BGL_n(\C)}  @>>> \MUL{B\glnk}
\end{CD}$$

The map on the left is an isomorphism because the variety of Grassmannians, ie $BGL_n(\C)$, has an algebraic cell decomposition, so $CH^*(BGL_n(\C))=H^*(BGL_n(\C),\Z)$ (cf \cite{fulton}), while on the other hand $H^*(BGL_n(\C),\Z)=\MU{BGL_n(\C)}$ because the cohomology ring has no torsion. This would be enough to conclude that $cl^0$ is surjective and hence an isomorphism if we could only be sure that the diagram is commutative. However this is not clear: the horizontal maps are essentially coming from the two different maps $B\glnk\to BGL_n(\C)$ constructed in \ref{blank: brauer} and \ref{morebrauer}, and it would require some work to prove that they coincide, if they do at all.

Instead we take a different route, which does not even use the computation {\em a priori} of $\MUL{B\glnk}$. It has the virtue of involving an interesting lemma, a souped-up version of the results of Ravenel-Wilson-Yagita \cite{bpfrommorava} and Hopkins-Kuhn-Ravenel \cite{moravaandgenchar}:

\begin{lem*} \label{lem: euler}
Let $G$ be a finite group such that $K(m)^*(BG)$ is generated by transferred Euler classes as a $K(m)^*$-module, for all $m$. Then $G\wr  S_n$ has the same property, for any $n\ge 0$. 
\end{lem*}

\begin{proof}
We shall make use of the fact that if $A$ is a subgroup of the finite group $H$ with index prime to $l$, then the map $K(m)^*(BH)\to K(m)^*(BA)$ is injective whereas the transfer is surjective (where, as always, we use Morava K-theory at the prime $l$). This is proved in \cite{moravaandgenchar}. We use the same type of induction as before.
\begin{enumerate}

\item{$n=n'l$.} In this case the subgroup $(G\wr  \Z/l)\wr 
S_{n'}$ has index prime to $l$, and the corestriction (transfer) map is surjective, so it is enough to prove the result for this group. But $G\wr  \Z/l$ satisfies the same property as $G$ by \cite{moravaandgenchar}, and we only have to apply the induction hypothesis.

\item{$n$ is prime to $l$.} The subgroup $G\times (G\wr  S_{n-1})$
has index prime to $l$, so again it is enough to prove the result for this subgroup. The Kunneth formula for Morava K-theories gives
$$K(m)^*(B(G\times G \wr  S_{n-1}))=K(m)^*(BG)\hat\otimes_{K^*(m)} K(m)^*(B(G \wr  S_{n-1}))$$ from which the result follows at once.

\end{enumerate}
This completes the induction step.\end{proof}

\begin{coro*}[1]
Let $G$ be as above. Then $\MUL{BG}$ is generated by transfered Euler classes, and so is $\MUL{B(G\wr  S_n)}$.

\end{coro*}

\begin{proof}
By corollary 2.2.1 in \cite{bpfrommorava}, $BP^*(BG)$ is generated by transfers of Euler classes as a $BP^*$-module, so $\BPL{BG}$ is generated by such classes as an abelian group, and the result follows for $G$. But by the lemma $G\wr  S_n$ can replace $G$.
\end{proof}

\begin{coro*}[2] \label{coro: cyclesurj}
For $G$ as above (for example, $G$ abelian), the two maps
$$CH^*(BG)/l\to \MUL{BG}$$ and $$CH^*(B(G\wr  S_n))/l\to \MUL{B(G\wr  S_n)}$$ are surjective.

\end{coro*}

\begin{rmk*} In this section ``generated'' means ``topologically generated'', and ``surjective'' means ``with a dense image''. The reader will check that this is of no consequence in what follows.
\end{rmk*}

\end{blank}

\begin{blank}{\em Conclusion.} We arrive finally at the promised isomorphism:
\begin{thm*} \label{thm: iso}
Let $k$ be a finite field of characteristic $p$ and let $l\ne p$ be a prime number.
Then the  map $$CH^*(B\glnk)/l \to \MUC{B\glnk}{\Z/l}$$ is an isomorphism.
\end{thm*}
\end{blank}

\begin{proof}
Recall that $\glnk$ has a subgroup $C^m$ (product of $m$ cyclic groups) which is normalized by some subgroup $N$, such that the restriction map induces an isomorphism $CH^*(B\glnk)/l\to (CH^*(BC^m)/l)^N$, the subring of fixed elements.
Now consider the following commutative diagram:

$$
\begin{CD}
CH^*(B\glnk)/l @>i^*_0>> CH^*(B(C\wr  S_m))/l @>j^*_0>> CH^*(BC^m)/l \\
@VVV @VclVV @V\phi V\approx V \\
\MUC{B\glnk}{\Z/l} @>i^*_1>> \MUC{B(C\wr  S_m)}{\Z/l} @>j^*_1>> \MUC{BC^m}{\Z/l}
\end{CD}
$$

 We have just indicated the image of the composition $j^*_0\circ i^*_0$.

The maps $i^*_0$ and $i^*_1$ are injective for index reasons, and $j^*_0$ is known to be injective. By (\ref{coro: cyclesurj}, corollary (2)), the map $cl$ is surjective, and it follows that $j^*_1$ is injective.

Also, because the isomorphism $\phi$ above ($C^m$ being abelian) is natural with respect to maps induced by homomorphisms of groups, it induces also an isomorphism between the subrings of fixed elements under the action of $N$. This describes the image of $j^*_1\circ i^*_1$ (namely, all of the fixed subring), and concludes the proof.\end{proof}

\section{A word on ordinary cohomology}

As pointed out in the introduction, the computation of the cobordism ring of a space $X$, or rather the simplified ring $\MUL{X}$, will give some information about the ordinary cohomology ring $H^*(X,\F_l)$. In the context of Chevalley groups as in (\ref{defchev}), it is easy to see that the map
$$\MUL{BG(K)}\longrightarrow H^*(BG(K),\F_l)$$ is injective, if $K$ is a finite field containing the $l$-th roots of unity. Indeed, for tori this is well-known (cf (\ref{intro})), and the general case follows from (\ref{injecttorus}) and (\ref{blank: bigones}). So we know that $H^*(BG(K),\F_l)$ contains a polynomial ring, which is also the image of $H^*(B\gc,\F_l)$ under the evident map (cf (\ref{blank: brauer})); incidentally this map is injective. It is reasonable to conjecture that this map becomes surjective after dividing each ring by its radical (ie the image of the map contains everything except some nilpotent elements). In fact we can expect the cohomology ring to be obtained from the cobordism ring by tensoring with an exterior algebra: this is the case for abelian groups and for $GL_n$, by \cite{quillengln} and the present article.

These cohomology rings have been investigated, see \cite{chevclass} and \cite{chevexcept}. Using more or less case-by-case considerations, we can see that the cohomology ring is in most cases the tensor product of some polynomial ring and some exterior algebra. Moreover the polynomial part is abstractly isomorphic to the cobordism ring (as graded rings, ie we have the same number of generators in the same degrees). So for most choices of $G$, $n$, $p$ and $l$ the conjecture above will be true. However it would still be interesting to have a direct, neat proof of this. The geometry beyond this problem involves comparing usual complex cobordism and ``cobordisms with singularities'' -- ordinary cohomology being the example in view.

\begin{acknowledgements}
This work was carried out at the University of Cambridge during the first year of my PhD. It is a great pleasure to thank my supervisor Burt Totaro for his guidance.
I am also indebted to D. Ravenel and S. Wilson for their help.
\end{acknowledgements}

\bibliography{myrefs}
\bibliographystyle{siam}

\end{document}